\input amstex
\documentstyle{amsppt}
\NoBlackBoxes

\hsize 6.2truein \vsize 8.9truein

\topmatter
\title Volume Inequalities for Isotropic Measures\endtitle
\author Erwin Lutwak, Deane Yang and Gaoyong Zhang \endauthor

       \leftheadtext{Volume Inequalities for Isotropic Measures }
       \rightheadtext{Lutwak, Yang, and Zhang}
       \author Erwin Lutwak, Deane Yang, and Gaoyong Zhang
       \endauthor
          \thanks Research supported, in part, by NSF Grant
          DMS--0104363 and DMS--0405707
          \endthanks

      \abstract A direct approach to Ball's {\it simplex inequality} is presented.
This approach, which does not use the Brascamp-Lieb inequality, also
gives Barthe's characterization of the simplex for Ball's inequality
and extends it from discrete to arbitrary measures. It also yields
the dual inequality, along with equality conditions, and it does
both for arbitrary measures.
    \endabstract
          \subjclass MR 52A40\endsubjclass
       \affil
        Department of Mathematics\\
        Polytechnic University\\
        Brooklyn, NY 11201\\
        \endaffil

\endtopmatter

\document

       \define\sn{S^{n-1}}
       \define\sm{S^{n}}

        \define\is{\int _{\sn}}

        \define\vv11{\tilde V_{-1}}

        \define\supp{\operatorname{supp}}
        
        \define\rp{\rangle_{\scriptscriptstyle\negthinspace\kern-.5pt p}}

        \define\rn{{\Bbb R}^n}
        \define\rmm{ {\Bbb R}^{n+1} }
        \define\rbo{\Bbb R}

\define\intt{\operatorname{int}}

\define\conv{\operatorname{conv}}

\define\inte{\operatorname{int}}

A non-negative Borel measure $Z$ on the unit sphere $S^{n-1}$ is
said to be isotropic if, when viewed as a mass distribution on
$\sn$, it has the same moment of inertia about all lines through the
origin. The convex hull of the support of $Z$ is a convex body
$Z_\infty$ in $\Bbb R^n$. Ball used his elegant reformulation of the
Brascamp-Lieb inequality to obtain the sharp upper bound for the
volume of the polar $Z_\infty^*$ of the body $Z_\infty$. Ball showed
how this could be combined with John's characterization of the
ellipsoid of maximal volume contained in a convex body, to prove his
celebrated reverse isoperimetric inequality \cite{2}. Recently,
Barthe \cite{5} found a beautiful new approach to establishing the
Brascamp-Lieb inequality. Barthe's approach has the added advantage
that it yields the equality conditions for Ball's reformulation of
the Brascamp-Lieb inequality. Armed with these equality conditions,
Barthe was able to show that, when $Z$ is discrete, the simplex is
the unique extremal for Ball's inequality.

In this paper we give a direct proof that for and only for the
simplex is the sharp upper bound for the volume of $Z_\infty^*$
attained. The advantage of our approach is that it allows us to
establish Ball's simplex inequality, along with Barthe's equality
conditions, for measures which are not necessarily discrete. Another
advantage of our approach is that it allows us to establish the
analogous results for the body $Z_\infty$.

The equality conditions of the Brascamp-Lieb inequality are
complicated. Our approach is based on a continuous version of the
Ball-Barthe inequality that has simple and easily stated equality
conditions.

We have attempted to write an article that is self-contained. We
have given very detailed proofs (but did not reprove the
Ball-Barthe inequality and refer the reader to \cite{28}).
Although the questions we address (here and in \cite{28}) arise
naturally within the Brunn-Minkowski-Firey theory (see, e\.g\.,
[7-10,13-27,29-32,34,35,37]) none of the machinery of the theory is used.
\medpagebreak

The ideas and techniques of Ball and Barthe play a critical role
in this paper. It would be impossible to overstate our reliance on
their work. \bigpagebreak

If $Z$ is a finite nonnegative Borel measure on the unit sphere
$S^{n-1}$ and is not concentrated on a closed hemisphere of $\sn$,
then define the convex body $Z_\infty$ in $\Bbb R^n$ as the convex
hull of the support of $Z$; i.e.,
$$
Z_\infty = \conv(\supp Z).\tag 1
$$
Note that the convex body $Z_\infty$ contains the origin in its
interior (since $\supp Z$ is not contained in a closed hemisphere of
$\sn$). The polar body $Z_\infty^*$ of $Z_\infty$ is given by
$$
Z_\infty^* =\{ x\in \Bbb R^n : x\cdot v \le 1\ \text{for all}\
v\in \supp Z\},\tag 2
$$
where $x\cdot v$ denotes the standard inner product of $x$ and $v$
in $\rn$. Note that since $Z_\infty$ is by definition contained in
the unit ball, $B$, it follows that $B\subseteq Z^*_\infty$.

In addition to its denoting absolute value, we shall use $|\cdot|$
to denote the standard Euclidean norm on $\rn$, on occasion the
absolute value of the determinant of an $n\times n$ matrix, and
often to denote $n$-dimensional volume.

A finite nonnegative Borel measure $Z$ on $S^{n-1}$ is said to be
{\it isotropic} if
$$
|x|^2 = \int_{S^{n-1}} | x\cdot u|^2 dZ(u),\tag 3
$$
for all $x\in\Bbb R^n$. Note that it is impossible for an
isotropic measure to be concentrated on a proper subspace of
$\rn$. The {\it centroid} of the measure $Z$ is defined as
$$
\frac1{Z(\sn)}\int_{S^{n-1}} u\, dZ(u).
$$
\bigpagebreak

The purpose of this paper is to establish the following two
theorems.

\proclaim{Theorem 1} If $Z$ is an isotropic measure on $\sn$
whose centroid is at the origin, then
$$
|Z_\infty^*| \le n^{n/2} (n+1)^{(n+1)/2}/n!,
$$
with equality if and only if $Z_\infty$ is a regular simplex
inscribed in $S^{n-1}$.
\endproclaim

The volume inequality was proved by Ball \cite{2}. For discrete
measures, the equality conditions were obtained by Barthe
\cite{5}.

The following dual to the inequality in Theorem 1
was anticipated by Ball \cite{2}
and established for discrete measures by Barthe \cite{4}.

\proclaim{Theorem 2} If $Z$ is an isotropic measure on $\sn$ whose
centroid is at the origin, then
$$
|Z_\infty| \ge (n+1)^{(n+1)/2} n^{-n/2}/ n!,
$$
with equality if and only if $Z_\infty$ is a regular simplex
inscribed in $S^{n-1}$.
\endproclaim

\bigpagebreak

Write $\{e_1,\ldots,e_n\}$ for the standard orthonormal basis of
$\rn$. If, in the definition (3) of an isotropic measure, we set
$x=e_j$ and sum over all $j$ we immediately see that the total mass
of an isotropic measure on $\sn$ is always $n$; i.e.,
$$
Z(\sn)=n.\tag 4
$$
We shall also make use of the simple observation that if $\supp
Z=\{u_1,\ldots,u_n\}\subset\sn$, then the vectors $u_1,\ldots,
u_n$ must be orthogonal. To see this first note that since $Z$ is
isotropic, from (3) we have for all $x\in\rn$,
$$
\sum_{i=1}^n c_i | x\cdot u_i|^2 =| x|^2,
$$
where $c_i=Z(\{u_i\})>0$. Taking $x=u_j$, gives
$$
   \sum_{i=1}^n c_i | u_j\cdot u_i|^2 =1.
$$
This shows that $c_j\le 1$. But from (4) we know $\sum_{i=1}^n c_i
= n$ and hence, $c_j=1$. Thus, $| u_j\cdot u_i|=0$ for $j\neq i$.
\bigpagebreak

If $\nu$ is a Borel measure on $\sn$, let $|f\negthinspace
:\negthinspace \nu|_p$ denote the standard $L_p$ norm of the
function $f$, with respect to $\nu$; i.e., for $1\le p<\infty$
$$
    |f\negthinspace :\negthinspace \nu|_p = \left(\is |f(u)|^p\, d\nu(u)    \right)^{1/p}.
$$

The {\it support function}  $h_K: \sn\to\rbo$ of a convex body $K$
in $\rn$ is defined by
$$
h_K (u) = \max \{\,u\cdot x : x\in K\}
$$
for $u \in\sn$. For each measure $Z$ on $\sn$ the convex body
$Z_\infty\subseteq B$ may be defined as the convex body whose
support function is given by
$$
h_{Z_\infty}(u) = \max \{\, u\cdot v : v\in \supp Z \},
$$
for $u\in\sn$. From (2) we see that
$$
\intt Z_\infty^* =\{ x\in \Bbb R^n : x\cdot v < 1\ \text{for all}\
v\in \supp Z\}.\tag 5
$$
\bigskip

\bigpagebreak

The following continuous version of the Ball-Barthe inequality was
proved in \cite{28} as a direct consequence of the H\"older
inequality.

\proclaim{Ball-Barthe lemma} If $\nu$ is an isotropic measure on
$S^{n-1}$ and $t : S^{n-1}\to(0,\infty)$ is continuous, then
$$
\det \int_{S^{n-1}} t(v)\, v\otimes v \, d\nu(v) \ge
\exp\bold{\{}\int_{S^{n-1}} \log t(v)\, d\nu(v)\bold{\}},
$$
with equality if and only if  $t(v_1)\cdots t(v_n)$ is constant for
linearly independent $v_1, \ldots, v_n\in \supp(\nu)$.
\endproclaim

\noindent Here, $v\otimes v:\rn\to\rn$ is the rank 1 linear operator
$x\mapsto (x\cdot v)v$.
\bigpagebreak

For $t\in L_1(Z)$, define $t^\circ \in\rn$ by
$$
t^\circ =\int_{S^{n-1}} u t(u)\, dZ(u).
$$

We shall make use of the following trivial fact.

\proclaim{Lemma 1} If $Z$ is an isotropic measure on $S^{n-1}$ and
$t\in L_2(Z)$, then
$$
|\,t^\circ |\le |\,t\negthinspace :\negthinspace Z|_2,
$$
with equality if and only if $t(u)=u\cdot t^\circ $ for almost all
$u\in\sn$ with respect to the measure $Z$.
\endproclaim
\demo{Proof} From the definition of $t^\circ $, the Cauchy-Schwarz
inequality, and finally using the fact that $Z$ is isotropic we have
$$
\align |\,t^\circ |^2=&\ t^\circ \negthinspace \cdot \,t^\circ \\
=&\is (t^\circ \cdot u)\,t(u)\,dZ(u)\\
\le&\is |\,t^\circ\negthinspace  \cdot u|\,|\,t(u)|\,dZ(u)\\
\le& |\,t\negthinspace :\negthinspace Z|_2\left[ \is |\,t^\circ\negthinspace  \cdot u|^2dZ(u)\right]^{\frac12}\\
=&|\,t\negthinspace :\negthinspace Z|_2\,|\,t^\circ |.
\endalign
$$\qed
\enddemo

We shall require the following simple fact:

\proclaim{Lemma 2} Suppose $Z$ is a measure on $\sn$ whose support
is not contained in a closed hemisphere of $\sn$. If
$t:\sn\to(0,\infty)$ is a continuous function such that
$|\,t\negthinspace :\negthinspace Z|_1=1$, then $t^\circ \in\inte
Z_\infty$.
\endproclaim
\demo{Proof} Suppose $v\in\sn$. Since $t>0$, and $\supp Z$ is not
contained in a closed hemisphere of $\sn$, we have
$$
v\cdot t^\circ =\is (v\cdot u)\, t(u)\,dZ(u) < |\,t\negthinspace
:\negthinspace Z|_1\,\max_{u\in\supp Z} v\cdot u=h_{Z_\infty}(v),
$$
for all $v\in\sn$, which gives the desired result that $t^\circ
\in\inte Z_\infty$.\qed
\enddemo

We shall make use of Lemma 2 only in the following form: If
$t:\sn\to(0,\infty)$ is continuous, then
$$
  \left.\is u t(u)\,dZ(u)\right/\is t(u)\,dZ(u)\ \in \inte Z_\infty.\tag 6
$$
\bigpagebreak

We shall need the function $s:\sn\to\sm\subset\rmm=\rn\times\rbo$,
defined by
$$
s(u)=(-\tsize\frac{\sqrt{n}}{ \sqrt{n+1} }\,u, \tsize\frac{1}{
\sqrt{n+1} }),\tag 7
$$
for each $u\in\sn$.

Following Ball's idea, we show how each isotropic measure $Z$ on
$\sn$, whose centroid is at the origin, induces an isotropic measure
${\bar Z}$ on $S^n$. Define the measure ${\bar Z}$ on $\sm$,
concentrated on the subsphere
$$
D=\{x\in\sm: x\cdot e_{n+1}=1/\sqrt{n+1}\},
$$
by
$$
\int_{S^n} f\, d{\bar Z} =  \int_D f\, d{\bar Z} =\frac{n+1}{n}
\int_{\sn} f\negthinspace\circ\negthinspace s\,dZ,\tag 8
$$
for each continuous $f:\sm\to\rbo$.

To see that ${\bar Z}$ is an isotropic measure on $S^n$, i.e.,
$$
\int_{S^n} |y\cdot v|^2 d{\bar Z}(v) = |y|^2,\tag 9
$$
for all $y\in \Bbb R^{n+1}$, note that from definition (8), the fact
that $Z$ is isotropic, the fact that the centroid of $Z$ is at the
origin, and (4), we have, for $y=(x,r)\in\Bbb
R^{n+1}=\rn\times\rbo$,
$$
 \align
  \int_{\sm} |y\cdot w|^2\,d\bar Z(w)
  &={\frac{n+1}{n}}\int_{\sn} |(x,r)\cdot
    (-\tsize\frac{\sqrt{n}}{\sqrt{n+1}}\,u,\tsize\frac{1}{\sqrt{n+1}})|^2\,dZ(u)\\
  &=\is |x\cdot u|^2\,dZ(u)\ -\
    \frac{2r}{\sqrt{n}}\, x\cdot\is u\,dZ(u)\ +\
    \frac{r^2}{n}\is dZ(u)\\
  &= |x|^2 + r^2\\
  &= |y|^2.
 \endalign
$$
Since $\bar Z$ is an isotropic measure on $\sm$, from (4) we have
$$
\bar Z(\sm)=n+1.\tag 10
$$
Observe that since $\bar Z$ is isotropic, $\supp\bar Z$ cannot be
contained in a subspace of $\rmm$.

From the fact that $Z$ has its centroid at the origin and the fact
that the plane of support of $\bar Z$ is orthogonal to $e_{n+1}$,
and passes through $(0,\ldots,0,1/\sqrt{n+1})$, it is easily seen
that the centroid of $\bar Z$ is $e_{n+1}/\sqrt{n+1}$; i.e.,
$$
\int_{\sm} w\,d\bar Z(w) = \sqrt{n+1}\,e_{n+1}.\tag 11
$$
Indeed, from definition (8), the fact that the centroid of $Z$ is at
the origin, and (4), we have for each $y=(x,r)\in\Bbb
R^{n+1}=\rn\times\rbo$,
$$
\align y\cdot \int_{\sm} w\,d\bar Z(w)
 &=\frac{n+1}{n}\int_{\sn} (x,r)\cdot
 (-\tsize\frac{\sqrt{n}}{\sqrt{n+1}}\,u,\tsize\frac{1}{\sqrt{n+1}})\,dZ(u)\\
 &=-{\tsize\frac{\sqrt{n+1}}{\sqrt{n}}}\, x\cdot \is  u\,dZ(u)
    + {\tsize\frac{\sqrt{n+1}}{n}}\,r\is dZ(u)\\
 &= \sqrt{n+1}\, r\\
 &= y\cdot (\sqrt{n+1}\,e_{n+1}).
\endalign
$$
\bigpagebreak

We now prove:

\proclaim{Theorem 1} If $Z$ is an isotropic measure on $S^{n-1}$
whose centroid is at the origin, then
$$
|Z_\infty^*| \le n^{n/2} (n+1)^{(n+1)/2}/n!,
$$
with equality if and only if $Z_\infty$ is a regular simplex
inscribed in $S^{n-1}$.
\endproclaim
\demo{Proof} Define the strictly increasing function $\phi :
(0,\infty) \to \Bbb R$ by
$$
\int_0^t e^{-s} ds = \frac1{\sqrt \pi} \int_{-\infty}^{\phi(t)}
e^{-s^2} ds.
$$
Note that $\phi'>0$ and that for all $t>0$, we have
$$
-t=-\log\sqrt{\pi} -\phi(t)^2 +\log\phi'(t).\tag 12
$$

Define the open cone $C\subset\rmm=\rn\times\rbo$ by
$$
C=\bigcup_{r>0} (\tsize\frac r{\sqrt n}\, \intt
Z_\infty^*)\times\{r\}.\tag 13
$$
Now $y=(x,r)\in C$ if and only if $\frac{\sqrt n}r\,x\in
\intt(Z_\infty^*)$, and from (5) we see that this is equivalent to $
x\cdot u < \frac r{\sqrt n}$ for all $u\in \supp Z$. Using (7) we
may rewrite this as
$$
y\in C\qquad\text{if and only if}\qquad y\cdot s(u)>0,\ \text{ for
all $u\in \supp Z$}.\tag 14
$$

Following Barthe's idea, define a transformation $T :C \to \Bbb
R^{n+1}$ by
$$
Ty=\int_{S^n} w\phi(y\cdot w)\, d{\bar Z}(w),\tag 15
$$
or, by (8), equivalently,
$$
Ty=\frac{n+1}n\is s(u)\phi(y\cdot s(u))\,dZ(u),\tag 16
$$
for each $y\in C$. From (14) we see that $y\cdot s(u)$ is in the
domain of $\phi$ for each $y\in C$ and each $u\in\supp Z$.

From (15) it follows that the differential of $T$ is given by
$$
dT(y)=\int_{S^n} \phi'(y\cdot w)\, w\otimes w\, d{\bar Z}(w),\tag 17
$$
for each $y\in C$. Thus, for each $z\in\rmm$,
$$
z\cdot dT(y)z=\int_{S^n} \phi'(y\cdot w) |w\cdot z|^2 d{\bar
Z}(w).
$$
Since $\phi'>0$ and $\bar Z$ is not concentrated on a proper
subspace of $\rmm$, we conclude that the matrix $dT(y)$ is positive
definite for each $y\in C$. Hence, a simple application of the mean
value theorem shows that $T:C\to\rmm$ is globally 1-1.

From (12), (10), (17) and the Ball-Barthe inequality, (15) and Lemma
1, and making the change of variables $z=Ty$, we have
$$
\align
 \int_C&\exp\left\{-\int_{S^n} y\cdot w\, d{\bar Z(w)}\right\} dy \\
 &=
 \int_C \exp\left\{-\int_{S^n} (\phi(y\cdot w)^2
       - \log\phi'(y\cdot w) + \log\sqrt{\pi})
       \,d{\bar Z(w)}\right\} dy \\
 &= \pi^{-\frac{n+1}2} \int_C
        \exp\left\{ -\int_{S^n} \phi(y\cdot w)^2 d{\bar Z(w)}  \right\}
        \exp\left\{ \int_{S^n} \log\phi'(y\cdot w) d{\bar Z(w)} \right\}
                       dy \\
 &\le \pi^{-\frac{n+1}2} \int_C\exp\left\{-\int_{S^n} \phi(y\cdot w)^2
       d{\bar Z(w)}\right\}|dT(y)|dy \\
 &\le \pi^{-\frac{n+1}2} \int_C\exp\left\{-|Ty|^2\right\}|dT(y)|dy \\
 &\le \pi^{-\frac{n+1}2} \int_{\Bbb R^{n+1}} e^{-|z|^2} dz \\
 &=   1.
\endalign
$$

On the other hand, from (11) and (13), we have
$$\align
 \int_C\exp\left\{-\int_{S^n} y\cdot w\, d{\bar Z(w)}\right\}dy
 &= \int_C\exp\left\{-y\cdot \int_{S^n} w\, d{\bar Z(w)}\right\}dy\\
 &= \int_C \exp\{-\sqrt{n+1}\,y\cdot e_{n+1}\}\,dy\\
 &=\int_0^\infty  \int_{\frac r{\sqrt n} \intt(Z_\infty^*)} e^{-\sqrt{n+1}\, r}dxdr \\
 &=|Z_\infty^*| \int_0^\infty \left(\frac r{\sqrt n}\right)^n e^{-\sqrt{n+1}\, r}dr \\
 &=|Z_\infty^*|\, n^{-n/2}\, n!\, (n+1)^{-(n+1)/2},
\endalign
$$
which gives the desired inequality. \medpagebreak

Suppose there is equality in the inequality of our theorem. Since
${\bar Z}$ is not concentrated on a proper subspace of $\rmm$,
there are linearly independent unit vectors $w_1,\ldots,w_{n+1}
\in \supp {\bar Z}$. Assume there is a different unit vector
$w_0\in\supp {\bar Z}$. Write $w_0=\lambda_1 w_1
+\cdots+\lambda_{n+1} w_{n+1}$. At least one coefficient, say
$\lambda_1$, is not zero. Since $w_0, w_2, \ldots, w_{n+1}$ are
linearly independent, the equality conditions of the Ball-Barthe
inequality imply that
$$
 \phi'(y\cdot w_1) \cdots \phi'(y\cdot w_{n+1})
 =
 \phi'(y\cdot w_0) \phi'(y\cdot w_2)\cdots \phi'(y\cdot w_{n+1}),
$$
for all $y\in C$. But $\phi'>0$, and hence we have
$$\phi'(y\cdot w_1)= \phi'(y\cdot w_0),$$
for all $y\in C$. Differentiating both sides with respect to $y$
shows that
$$
\phi''(y\cdot w_1)\,w_1 = \phi''(y\cdot w_0)\,w_0,
$$
for all $y\in C$. Since there exists $y \in C$ such that
$\phi''(y\cdot w_1)\neq 0$ it follows that $w_0 = \pm w_1$. But $Z$
is supported inside an open hemisphere of $S^n$, so $w_0 = w_1$.

Hence equality in our inequality implies that $\supp {\bar
Z}=\{w_1,\ldots, w_{n+1}\}$. Since ${\bar Z}$ is isotropic,
$w_1,\ldots, w_{n+1}$ are orthogonal. But $w_i\perp w_j$ implies
$u_i\cdot u_j = 1/n$, for $i\neq j$, and thus $\supp Z$ consists of
the vertices of a regular simplex inscribed in $S^{n-1}$. \qed
\enddemo

When $Z$ is a discrete measure, the inequality of Theorem 1 was
proved by Ball \cite{2}. For discrete measures, the equality
conditions of Theorem 1 were obtained by Barthe \cite{5}.

The inequality of the following theorem was anticipated by Ball
\cite{2} and, for discrete measures, established by Barthe \cite{4}.

\proclaim{Theorem 2} If $Z$ is an isotropic measure on $S^{n-1}$
whose centroid is at the origin, then
$$
|Z_\infty| \ge (n+1)^{(n+1)/2} n^{-n/2}/ n!,
$$
with equality if and only if $Z_\infty$ is a regular simplex
inscribed in $S^{n-1}$.
\endproclaim
\demo{Proof} Define the strictly increasing function $\phi : \Bbb
R \to (0,\infty)$ by
$$
\int_0^{\phi(t)} e^{-s} ds = \frac1{\sqrt \pi} \int_{-\infty}^{t}
e^{-s^2} ds.
$$
Note that $\phi'>0$, and that for all $t\in\rbo$, we have
$$
 t^2 = \phi(t) -\log\phi'(t)- \log\sqrt{\pi}.
\tag 18
$$

Define a transformation $T : \Bbb R^{n+1} \to \Bbb R^{n+1}$ by
$$
Ty=\int_{S^n} w\phi(y\cdot w)\, d{\bar Z}(w), \tag 19
$$
or, by (8), equivalently by
$$
Ty=\frac{n+1}n \is s(u)\phi(y\cdot s(u))\,dZ(u). \tag 20
$$
Define the open cone $C\subset\rmm=\rn\times\rbo$ by
$$
C=\bigcup_{r>0}  (-r\sqrt{n} \intt Z_\infty)\times\{r\}. \tag 21
$$
Note that if $z\in\rmm$ is such that $z\cdot e_{n+1} > 0$, then
$$
z\in C\qquad\text{if and only if}\qquad z|_{\rn}/ z\cdot e_{n+1} \in
-\sqrt{n}Z_\infty, \tag 22
$$
where $z|_{\rn}$ denotes the orthogonal projection of $z$ onto
$\rn$. We now show that $T(\rmm)\subset C$. To see this, note that
from (20) and definition (7) it follows that
$$
\left.\phantom{-\tsize\frac1{\sqrt{n}}\,}
Ty\right|_{\rn}=\,-\frac{\sqrt{n+1}}{\sqrt{n}} \int_{\sn} u
\phi(y\cdot s(u))\, d{Z(u)}, \tag 23
$$
while the $(n+1)$-st component of $Ty$ is given by
$$
 Ty\cdot e_{n+1}=
 \frac{\sqrt{n+1}}n \int_{\sn} \phi(y\cdot s(u))\,d{Z(u)},
 \tag 24
$$
or equivalently by
$$
 Ty\cdot e_{n+1}
 =
 \frac1{\sqrt{n+1}}\int_{\sm}\phi(y\cdot w)\,d{\bar Z}(w).
 \tag 25
$$
Since $\phi>0$, it follows from (23) and (24), together with (6),
that
$$
\left.\left.-\tsize\frac1{\sqrt{n}}\,Ty\right|_{\rn}\,\right/\,
Ty\cdot e_{n+1}\ \in \intt Z_\infty,
$$
and hence from (22) we have
$$
T(\Bbb R^{n+1})\subset C.\tag 26
$$

From (19) it follows that $dT$, the differential of $T$, is given by
$$
dT(y)=\int_{S^n}  w\otimes w\, \phi'(y\cdot w)\,d{\bar Z}(w),
 \tag 27
$$
for each $y\in\Bbb R^{n+1}$. Since $\phi'>0$ and $\bar Z$ is not
concentrated on a proper subspace of $\rmm$, we conclude from (27)
that the matrix $dT(y)$ is positive definite for each $y\in{\Bbb
R^{n+1}}$. Hence, an application of the mean value theorem shows
that the transformation $T:\Bbb R^{n+1} \to C$ is globally 1-1.

From (9), (18), (10), (27) together with the Ball-Barthe inequality
(25), making the change of variables $z=Ty$ and using (26), and
(21), we have
$$\align
\pi^\frac{n+1}2 &=\int_{\Bbb R^{n+1}} e^{-|y|^2}dy \\
 &=\int_{\Bbb R^{n+1}}\exp\left\{-\int_{S^n}|y\cdot w|^2 d{\bar Z(w)}\right\}dy \\
 &=\int_{\Bbb R^{n+1}}\exp\left\{-\int_{S^n} (\phi(y\cdot w)
    - \log\phi'(y\cdot w) - \log\sqrt{\pi})\,d{\bar Z(w)}\right\}dy \\
 &= \pi^\frac{n+1}2
    \int_{\Bbb R^{n+1}}
      \exp\left\{-\int_{S^n}\phi(y\cdot w)d{\bar Z(w)}\right\}
      \exp\left\{\int_{S^n}\log \phi'(y\cdot w) d{\bar Z(w)}\right\}
    dy\\
 &\le \pi^\frac{n+1}2
    \int_{\Bbb R^{n+1}}
       \exp\left\{-\int_{S^n} \phi(y\cdot w) d{\bar Z(w)}\right\}
    |dT(y)|\, dy \\
 &= \pi^\frac{n+1}2
    \int_{\Bbb R^{n+1}}
       \exp\{-\sqrt{n+1}\,\, Ty\cdot e_{n+1}\}
     |dT(y)|\, dy \\
 &\le \pi^\frac{n+1}2
    \int_{C} \exp\{-\sqrt{n+1}\,\, z\cdot e_{n+1}\}\, dz \\
 &\le \pi^\frac{n+1}2
    \int_0^\infty
    \int_{(-r\sqrt{n} Z_\infty\!)} e^{-\sqrt{n+1}\,\,r}\, dx\,dr \\
 &= \pi^\frac{n+1}2  n^{\frac{n}{2}}|Z_\infty|
    \int_0^\infty r^n e^{-\sqrt{n+1}\,\,r}
    \,dr \\
 &=\pi^\frac{n+1}2  n^{\frac{n}{2}}\,n!\,(n+1)^{-(n+1)/2}\,|Z_\infty|.
\endalign$$

Suppose there is equality in the inequality of our theorem. Since
${\bar Z}$ is not concentrated on a proper subspace of $\rmm$, there
are linearly independent unit vectors $w_1,\ldots,w_{n+1} \in \supp
{\bar Z}$. Assume $w_0$ is a different unit vector in $\supp {\bar
Z}$. Write $w_0=\lambda_1 w_1 +\cdots+\lambda_{n+1} w_{n+1}$. At
least one coefficient, say $\lambda_1$, is not zero. Since $w_0,
w_2, \ldots, w_{n+1}$ are linearly independent, the equality
conditions of the Ball-Barthe inequality imply that
$$
 \phi'(y\cdot w_1) \cdots \phi'(y\cdot w_{n+1})
 =
 \phi'(y\cdot w_0) \phi'(y\cdot w_2)\cdots \phi'(y\cdot w_{n+1}),
$$
for all $y \in \Bbb R^{n+1}$. Since $\phi'>0$, we have
$$\phi'(y\cdot w_1)= \phi'(y\cdot w_0),$$
for all $y \in \Bbb R^{n+1}$.

Since $\phi'$ is not constant, there are $c_1$ and $c_0$ so that
$\phi'(c_1)\neq \phi'(c_0)$. Since $\bar{Z}$ is supported in an open
hemisphere of $S^n$, the unit vectors $w_0$ and $w_1$ are not
parallel, and therefore there exists $y \in \rbo^{n+1}$ such that
$$y\cdot w_1=c_1,\qquad\text{and} \qquad y\cdot w_0=c_0,
$$
producing a contradiction.

Thus equality in the inequality implies that $\supp {\bar
Z}=\{w_1,\ldots, w_{n+1}\}$. Since ${\bar Z}$ is isotropic,
$w_1,\ldots, w_{n+1}$ are orthogonal, and thus, $\supp Z$ consists
of the vertices of a regular simplex inscribed in $S^{n-1}$. \qed
\enddemo

After a copy of this paper was communicated to him, Barthe \cite{6}
showed how these inequalities could also be obtained from a new
\lq\lq continuous\rq\rq\ version of the Brascamp-Lieb inequality.

\end